\documentclass[a4paper]{article}
\usepackage{latexsym}
\usepackage{amssymb}

\newtheorem{Theorem} {Theorem} [section]
\newtheorem{Proposition} [Theorem] {Proposition}
\newtheorem{Lemma} [Theorem] {Lemma}
\newcommand{\Proof}{ \noindent{\bf Proof.}\quad }
\newcommand{\qed}{\hfill$\Box$\medskip}
\newcommand{\Nn}{{\mathbb N}}
\newcommand{\Rr}{{\mathbb R}}
\newcommand{\Zz}{{\mathbb Z}}
\newcommand{\isom}{\,\simeq\,}

\newcommand{\maysplit}{\discretionary{}{}{}}
\newcommand{\presetctr}[1]{\renewcommand{\theTheorem}{\ref{#1}}}
\newcommand{\resetctr}{\addtocounter{Theorem}{-1}}

\title{Notes on simplicial rook graphs}
\author{Andries E. Brouwer, Sebastian M. Cioab\u{a}$^{*1}$,\\
Willem H. Haemers$^\dagger$ \& Jason R. Vermette$^*$\vspace{10pt}\\
\small $^*$Department of Mathematical Sciences,\vspace{-3pt}\\
\small University of Delaware, Newark DE 19716-2553, USA\vspace{4pt}\\
\small $^\dagger$Department of Econometrics and Operations Research,\vspace{-3pt}\\
\small Tilburg University, Tilburg, The Netherlands\vspace{4pt}\\
\small {\tt aeb@cwi.nl, cioaba@math.udel.edu,}\vspace{-3pt}\\
\small {\tt haemers@uvt.nl, vermette@math.udel.edu}
}

\begin{document}
\maketitle
\footnotetext[1]{
Research partially supported by the National Security Agency grant
H98230-13-1-0267.}
\vspace{-10pt}\begin{center}
{\it To Chris Godsil, on the occasion of his 65th birthday}
\end{center}

\begin{abstract}
The simplicial rook graph ${\rm SR}(m,n)$ is the graph of which the
vertices are the sequences of nonnegative integers of length $m$
summing to $n$, where two such sequences are adjacent when they differ
in precisely two places.
We show that ${\rm SR}(m,n)$ has integral eigenvalues,
and smallest eigenvalue $s = \max (-n, -{m \choose 2})$,
and that this graph has a large part of its spectrum in common with
the Johnson graph $J(m+n-1,n)$.
We determine the automorphism group and several other properties.
\end{abstract}

\section{Introduction}
Let $\Nn$ be the set of nonnegative integers, and let $m,n \in \Nn$.
The {\em simplicial rook graph} ${\rm SR}(m,n)$ is the graph $\Gamma$
obtained by taking as vertices the vectors in $\Nn^m$ with coordinate sum $n$,
and letting two vertices be adjacent when they differ in precisely
two coordinate positions.
Then $\Gamma$ has $v = {n+m-1 \choose n}$ vertices,
and is regular of valency $k = n(m-1)$.

Small cases:
For $n>0$, $m=0$, the graph is $K_0$.
For $n=0$ or $m=1$, the graph is $K_1$.
For $n=1$, the graph is $K_m$.
For $m=2$, the graph is $K_{n+1}$.
For $n=2$, the graph is the triangular graph $T(m+1)$.

The graph ${\rm SR}(m,n)$ was studied in detail by Martin \& Wagner \cite{MW}.
Here we settle their main conjecture, and show that this graph
has integral spectrum. Some results:

{\presetctr{integr_thm}
\begin{Theorem}
The graph ${\rm SR}(m,n)$ has integral spectrum.
\end{Theorem}
\resetctr}

{\presetctr{smallestev_prop}
\begin{Proposition}
If $\,{\rm SR}(m,n)$ has at least one vertex,
it has smallest eigenvalue $\max ( -n, -{m \choose 2})$.
\end{Proposition}
\resetctr}

The Johnson graph $J(m,n)$ (with ${m \choose n}$ vertices,
valency $n(m-n)$ and eigenvalues $(n-i)(m-n-i)-i$ with multiplicity
${m \choose i} - {m \choose i-1}$ $(0 \le i \le n)$, cf.~\cite{BCN})
is an induced subgraph of $\Gamma$ (on the set of vertices in $\{0,1\}^m$).

The graphs ${\rm SR}(m,n)$ and $J(m+n-1,n)$ both have
${m+n-1 \choose n}$ vertices and valency $n(m-1)$.
These graphs resemble each other and have a large part of
their spectrum in common.

{\presetctr{commonev}
\begin{Proposition}
The graphs ${\rm SR}(m,n)$ and $J(m+n-1,n)$ have equitable partitions
with the same quotient matrix $E$, where $E$ has eigenvalues
$(n-i)(m-i)-n$ with multiplicity ${m \choose i}$ for $0 \le i < n$,
and multiplicity ${m \choose n} - 1$ for $i = n$.
In particular, the spectrum of $E$ is part of the spectrum of
${\rm SR}(m,n)$ and $J(m+n-1,n)$.
\end{Proposition}
\resetctr}

Clearly, ${\rm Sym}(m)$ acts as a group of automorphisms,
permuting the coordinate positions. Except for small cases,
this is the full group.

{\presetctr{automgp}
\begin{Proposition}
Let $m,n > 2$. If $n > 3$, then ${\rm Aut} (\Gamma) \isom {\rm Sym}(m)$.
If $n = 3$, then ${\rm Aut} (\Gamma) \isom {\rm Sym}(m).2$.
\end{Proposition}
\resetctr}

We determine some miscellaneous properties.

{\presetctr{diameter}
\begin{Proposition}
Let $m > 0$. Then the diameter of ${\rm SR}(m,n)$ is $\min (m-1,n)$.
\end{Proposition}
\resetctr}

{\presetctr{cliquesz}
\begin{Proposition}
Let $m > 1$ and $n > 0$.
Then the size of the largest clique in ${\rm SR}(m,n)$ is
$\max(m,n+1)$.
\end{Proposition}
\resetctr}

We also give the independence number in case $n=3$,
and the complete spectrum for a few small $n$.
In general, the graphs ${\rm SR}(m,n)$ are not determined
by their spectrum.

{\presetctr{cospectral}
\begin{Proposition}
The graph ${\rm SR}(m,n)$ is not determined by its spectrum
when (a)~$m=4$ and $n \ge 3$, or (b)~$n = 3$ and $m \ge 4$.
\end{Proposition}
\resetctr}

Finally, we study the structure of the eigenspace for the eigenvalue $-n$.

\section{Integrality of the eigenvalues}
We start with the main result, proved by a somewhat tricky induction.

\begin{Theorem} \label{integr_thm}
All eigenvalues of ${\rm SR}(m,n)$ are integers.
\end{Theorem}
\Proof
Let $\Gamma$ be the graph ${\rm SR}(m,n)$, and let $X$ be its vertex set.
The adjacency matrix $A$ of $\Gamma$ acts as a linear operator on $\Rr^X$
(sending each vertex to the sum of its neighbors).
By induction we construct a series of subspaces
$0=U_0 \subseteq U_1 \subseteq ... \subseteq U_t = \Rr^X$
and find integers $c_i$, such that $(A-c_iI)U_i \subseteq U_{i-1}$
$(1 \le i \le t)$. Then $p(A) := \prod_i (A-c_iI)$ vanishes identically,
and all eigenvalues of $A$ are among the integers $c_i$.

For $j \ne k$, let $A_{jk}$ be the matrix that describes adjacency between
vertices that differ only in the $j$- and $k$-coordinates.
Then $A = \sum A_{jk}$. If $(A_{jk} - c_{jk}I)u \in U$ for all $j,k$
then $(A-cI)u \in U$ for $c = \sum c_{jk}$.

A basis for $\Rr^X$ is given by the vectors $e_x$ ($x \in X$)
that have $y$-coordinate 0 for $y \ne x$, and $x$-coordinate 1.
For $S \subseteq X$, let $e_S := \sum_{x \in S} e_x$,
so that $e_X$ is the all-1 vector. Since $(A-cI) e_X = 0$ for
$c = (m-1)n$, we can put $U_1 = \langle e_X \rangle$.

For partitions $\Pi$ of the set of coordinate positions $\{1,...,m\}$
and integral vectors $z$ indexed by $\Pi$ that sum to $n$, let
$S_{\Pi,z}$ be the set of all $u \in X$ with
$\sum_{i \in \pi} u_i = z_\pi$ for all $\pi \in \Pi$.
If $\Pi$ is a partition into singletons, then $|S_{\Pi,z}| = 1$.

For a vector $y$ indexed by a partition $\Pi$, let $\tilde{y}$ be the
sequence of pairs $(y_\pi,|\pi|)$ ($\pi \in \Pi$) sorted
lexicographically: with the $y_\pi$ in nondecreasing order,
and for given $y_\pi$ with the $|\pi|$ in nondecreasing order.

Order pairs $(\Pi,z)$ by $(\Sigma,y) < (\Pi,z)$ when $|\Sigma| < |\Pi|$,
or when $|\Sigma| = |\Pi|$ and $\tilde{y} \ne \tilde{z}$
and in the first place $j$ where $\tilde{y}$ and $\tilde{z}$ differ,
the pair $\tilde{y}_j$ is lexicographically smaller than the pair
$\tilde{z}_j$.

We use induction to show for $S = S_{\Pi,z}$ and suitable $c$
that the image $(A-cI) e_S$ lies in the subspace $U$ spanned
by $e_T$ for $T = S_{\Sigma,y}$, where $(\Sigma,y) < (\Pi,z)$.

Note that the sets $S = S_{\Pi,z}$ induce regular subgraphs of $\Gamma$.
Indeed, the induced subgraph is a copy of the Cartesian product
$\prod_\pi {\rm SR}(|\pi|,z_\pi)$.
The image $(A-cI)e_S$ can be viewed as a multiset where the $x \in X$
occur with certain multiplicities. The fact that $S$ induces a regular
subgraph means that we can adjust $c$ to give all $x \in S$ any desired
given multiplicity, while the multiplicity of $x \not\in S$ does not
depend on $c$.

If $j,k$ belong to the same part of $\Pi$, then $A_{jk} e_S$ only
contains points of $S$, and can be ignored.
So, let $j \in \pi$, $k \in \rho$, where $\pi,\rho \in \Pi$, $\pi \ne \rho$,
and consider $A_{jk} e_S$. Abbreviate $\pi \cup \{k\}$ with $\pi+k$
and $\pi \setminus \{j\}$ with $\pi-j$.

The image $(A_{jk}-cI)e_S$ equals $S_1-S_2$, where
$S_1$ is the sum of all $e_T$, with $T = S_{\Sigma,y}$ and
$\Sigma = (\Pi\setminus\{\pi,\rho\}) \cup \{\pi-j,\rho+j\}$
(omitting $\pi-j$ if it is empty)
and $y$ agrees with $z$ except that $y_{\pi-j} \le z_\pi$
and $y_{\rho+j} \ge z_\rho$ (of course $y_{\pi-j}+y_{\rho+j} =
z_\pi+z_\rho$), and $S_2$ is the sum of all $e_T$,
with $T = S_{\Sigma,y}$ and
$\Sigma = (\Pi\setminus\{\pi,\rho\}) \cup \{\pi+k,\rho-k\}$
and $y$ agrees with $z$ except that $y_{\pi+k} < z_\pi$
and $y_{\rho-k} > z_\rho$.

(Let $u$ be a $(j,k)$-neighbor of $s \in S$.
Since $\sum_{i \in \pi} s_i = z_\pi$, it follows that
$\sum_{i \in \pi-j} u_i = \sum_{i \in \pi-j} s_i \le z_\pi$,
so that $u$ is counted in $S_1$. Conversely, if $u$ is counted
in $S_1$, then we find a $(j,k)$-neighbor $s \in S$
by moving $u_j-s_j$ from position $j$ to position $k$
(if $u_j > s_j$) or moving $s_j-u_j$ from position $k$ to
position~$j$ (if $s_j > u_j$). The latter is impossible if
$u_k < s_j-u_j$, i.e. $\sum_{i \in \pi+k} u_i < z_\pi$,
and these cases are subtracted in $S_2$.)

We are done by induction. Indeed, for the pair $\{j,k\}$ we can
choose which of the two is called $j$, and we pick notation such that
$(z_\pi,|\pi|) \le (z_\rho,|\rho|)$ in lexicographic order.
Now in $S_1$ and $S_2$ only $(\Sigma,y)$ occur with
$(\Sigma,y) < (\Pi,z)$.
\qed

\section{The smallest eigenvalue}
We find the smallest eigenvalue of $\Gamma$ by observing
that $\Gamma$ is a halved graph of a bipartite graph $\Delta$.

Consider the bipartite graph $\Delta$ of which the vertices are
the vectors in $\Nn^m$ with coordinate sum at most $n$, where
two vertices are adjacent when one has coordinate sum $n$,
the other coordinate sum less than $n$, and both differ in precisely
one coordinate. Let $V$ be the set of vectors in $\Nn^m$ with
coordinate sum $n$. Two vectors $u,v$ in $V$ are adjacent in $\Gamma$
precisely when they have distance 2 in $\Delta$. If the adjacency matrix of
$\Delta$ is $\bigl({0 \atop N^\top}{N \atop 0}\bigr)$, with top and left
indexed by $V$, then for the adjacency matrix $A$ of $\Gamma$ we find
$A+nI = NN^\top$, so that $A+nI$ is positive semidefinite,
and the smallest eigenvalue of $A$ is not smaller than $-n$.

Together with the results of \cite{MW}, this proves that the smallest
eigenvalue of $A$ equals $\max (-n, -{m \choose 2})$.

\begin{Proposition} \label{smallestev_prop}
If \,${\rm SR}(m,n)$ has at least one vertex,
it has smallest eigenvalue $\max ( -n, -{m \choose 2})$.
\end{Proposition}
\Proof
Let $s$ be the smallest eigenvalue of $A$.
We just saw that $s \ge -n$.
Elkies \cite{Elkies} observed that $s \ge -{m \choose 2}$,
since $A$ is the sum of ${m \choose 2}$ matrices $A_{jk}$
that describe adjacency where only coordinates $j,k$ are changed.
Each $A_{jk}$ is the adjacency matrix of a graph that is a union of cliques,
and hence has smallest eigenvalue not smaller than $-1$.
Then $A = \sum A_{jk}$ has smallest eigenvalue not smaller than
$-{m \choose 2}$.

It is shown in \cite{MW} that the eigenvalue $-{m \choose 2}$
has multiplicity at least ${n-{m-1 \choose 2} \choose m-1}$,
and hence occurs with nonzero multiplicity if $n \ge {m \choose 2}$.
It is also shown in \cite{MW} that the multiplicity of the eigenvalue $-n$
is at least the number of permutations in ${\rm Sym}(m)$ with
precisely $n$ inversions, that is the number of words $w$ of length $n$
in this Coxeter group, and this is nonzero precisely when
$n \le {m \choose 2}$.
\qed

\begin{Proposition} \label{exactmult}
The eigenvalue $-{m \choose 2}$ has multiplicity precisely
${n-{m-1 \choose 2} \choose m-1}$.
\end{Proposition}
\Proof
For each vertex $u$, and $1 \le j < k \le m$, let $C_{jk}(u)$ be
the $(j,k)$-clique on $u$, that is the set of all vertices $v$
with $v_i = u_i$ for $i \ne j,k$.
An eigenvector $a = (a_u)$ for the eigenvalue $-{m \choose 2}$ must be
a common eigenvector of all $A_{jk}$ for the eigenvalue $-1$.
That means that $\sum_{v \in C} a_v = 0$ for each set $C = C_{jk}(u)$.

Order the vertices by $u > v$ when $u_d > v_d$
when $d = d_{uv}$ is the largest index where $u,v$ differ.
Suppose $u_i = s$ for some index $i$ and $s \le m-i-1$.
We can express $a_u$ in terms of $a_v$ for smaller $v$ with
$d_{uv} \ge m-s$ via $\sum_{v \in C} a_v = 0$, where $C = C_{i,m-s}(u)$.
Indeed, this equation will express $a_u$ in terms of $a_v$ where
$u_i + u_{m-s} = v_i + v_{m-s}$ and $v_j = u_j$ for $j \ne i, m-s$.
If $v_i > s$ this is OK since $v_{m-s} < u_{m-s}$.
If $t = v_i < s$, then by induction $a_v$ in its turn can be expressed
in terms of $a_w$ where $w$ is smaller and $d_{vw} \ge m-t > m-s$,
so that $w$ is smaller than $u$, and $d_{uw} > m-s$.

In this way we expressed $a_u$ when $u_i \le m-i-1$ for some $i$.
The free $a_u$ have $u_i \ge m-i$ for all $i$, and the vector $u'$ with
$u'_i = u_i - (m-i)$ is nonnegative and sums to $n - {m \choose 2}$.
There are
${n - {m \choose 2} + m-1 \choose m-1} = {n-{m-1 \choose 2} \choose m-1}$
such vectors, so this is an upper bound for the multiplicity.
But by \cite{MW} this is also a lower bound.
\qed

Thanks to a suggestion by Aart Blokhuis, we can also settle
the multiplicity of the eigenvalue $-n$.

\begin{Proposition}
The multiplicity of the eigenvalue $-n$ equals
the number of elements of ${\rm Sym}(m)$ with $n$ inversions,
that is, the coefficient of $t^n$ in the product
$\prod_{i=2}^m (1 + t + \cdots + t^{i-1})$.
\end{Proposition}
\Proof
As already noted, it is shown in \cite{MW} that the multiplicity
of the eigenvalue $-n$ is at least the number of permutations
in ${\rm Sym}(m)$ with precisely $n$ inversions.
The formula $\sum_w t^{l(w)} = \prod_{i=2}^m {t^i-1 \over t-1}$,
where $w$ runs over ${\rm Sym}(m)$ and $l(w)$ is the number of
inversions of $w$, is standard, cf.~\cite{Humphreys}, p.\,73.

Since $A+nI = NN^\top$, the multiplicity of the eigenvalue $-n$
is the nullity of $N$, and we need an upper bound for that.

\medskip
We first define a matrix $P$, and observe that
$N$ and $P$ have the same column space, and hence the same rank.
For $u,v \in \Nn^m$, write $u \preceq v$ when $u_i \le v_i$ for all $i$.
Let $P$ be the 0-1 matrix with the same row and column indices
(elements of $\Nn^m$ with sum $m$ and sum smaller than $m$, respectively)
where $P_{xy} = 1$ when $y \preceq x$.
Recall that $N$ is the 0-1 matrix with $N_{xy} = 1$ when $x$ and $y$
differ in precisely one coordinate position.
Let $M(y)$ denote column $y$ of the matrix $M$.

For $d = n - \sum y_i$ we find that
$$N(y) = \sum_{i=0}^{d-1} (-1)^i (d-i) \sum_{z \in W_i} P(y+z),$$
where $W_i$ is the set of vectors in $\{0,1\}^m$ with sum $i$.
Indeed, suppose that $x$ and $y$ differ in $j$ positions.
Then $j \le d$, and $N_{xy} = \delta_{1j}$, while the $x$-entry
of the right hand side is $\sum_{i=0}^d (-1)^i (d-i) {j \choose i} =
j \sum_{i=1}^j (-1)^{i-1} {j-1 \choose i-1} = j(1-1)^{j-1} = \delta_{1j}$.
We see that $N(y)$ and $dP(y)$ differ by a linear combination
of columns $P(y')$ where $\sum y'_i > \sum y_i$, and hence
that $N$ and $P$ have the same column space.

\medskip
Aart Blokhuis remarked that the coefficient of $t^n$ in the product
$\prod_{i=2}^m (1 + t + \cdots + t^{i-1})$ is precisely the number
of vertices $u$ satisfying $u_i < i$ for $1 \le i \le m$.
Thus, it suffices to show that the rows of $N$ (or $P$) indexed by
the remaining vertices are linearly independent.

Consider a linear dependence between the rows of $P$ indexed
by the remaining vertices, and let $P'$ be the submatrix of $P$
containing the rows that occur in this dependence.
Order vertices in reverse lexicographic order, so that $u$ is earlier
than $u'$ when $u_h < u'_h$ and $u_i = u'_i$ for $i > h$.
Let $x$ be the last row index of $P'$ (in this order).
Let $h$ be an index where the inequality $x_i < i$ is violated,
so that $x_h \ge h$.
Let $e_i$ be the element of $\Nn^m$ that has all coordinates 0
except for the $i$-coordinate, which is 1.
Let $z = x - he_h$.
Let $H = \{1,\ldots,h-1\}$. For $S \subseteq H$, let
$\chi(S)$ be the element of $\Nn^m$ that has $i$-coordinate 1
if $i \in S$, and 0 otherwise.

Consider the linear combination $p = \sum_S (-1)^{|S|} P'(z + \chi(S))$
of the columns of $P'$. We shall see that $p$ has $x$-entry 1
and all other entries equal to 0.
But that contradicts the existence of a linear dependence.

If $u$ is a row index of $P'$, and not $z \preceq u$, then $p_u = 0$.
If $z \preceq u$, and $z_i < u_i$ for some $i < h$, then the alternating sum
vanishes, and $p_u = 0$.
So, if $p_u \ne 0$, then $u$ agrees with $x$ in coordinates $i$,
$1 \le i \le h-1$.
For row $x$ only $S = \emptyset$ contributes, and $p_x = 1$.
Finally, if $u \ne x$ then $u_i \ge x_i$ for $i \ne h$ and
$\sum x_i = \sum u_i = n$ imply that $u_h < x_h$ and $u_i > x_i$
for some $i > h$, which is impossible, since $x$ is the reverse
lexicographically latest row index of $P'$.
\qed

These three propositions settle conjectures from \cite{MW}.

\section{An equitable partition}
A partition $\{X_1,\ldots,X_t\}$ of the vertex set $X$ of a graph $\Gamma$
is called {\em equitable} when for all $i,j$ the number $e_{ij}$ of
vertices in $X_j$ adjacent to a given vertex $x \in X_i$ does not depend on
the choice of $x \in X_i$.
In this case the matrix $E = (e_{ij})$ is called the {\em quotient matrix}
of the partition. All eigenvalues of $E$ are also eigenvalues
of $\Gamma$, realized by eigenvectors that are constant on the sets $X_i$.
There is a basis of $\Rr^X$ consisting of eigenvectors that either
are constant on all $X_i$, or sum to zero on all $X_i$.
The partition of $X$ into orbits of an automorphism group $G$ of $\Gamma$
is always equitable.

In this section we indicate an equitable partition of ${\rm SR}(m,n)$,
and in the next section a much finer one.

\medskip
Let $\Gamma$ be the graph ${\rm SR}(m,n)$ where $n > 0$, and let
$V_i$ be the set of vertices with precisely $i$ nonzero coordinates,
so that $|V_i| = {m \choose i} {n-1 \choose i-1}$.

\begin{Proposition}
The partition $\{V_1,\ldots,V_{\min(m,n)}\}$ is equitable.
Each $x \in V_i$ has $i(i-1)$ neighbors in $V_{i-1}$,
$(n-i)(m-i)$ neighbors in $V_{i+1}$, and all other neighbors in $V_i$.
The quotient matrix $E$ is tridiagonal, and has eigenvalues
$(m-i)(n-i) - n$ for $0 \le i \le \min(m,n)-1$.
\end{Proposition}
\Proof
The $e_{ij}$ are easily checked. It remains to find the eigenvalues.
Let $u$ be an eigenvector for the Johnson graph $J(m+n-1,n)$
for the eigenvalue $\theta = (m-i)(n-i)-n$.
Then $c_i u_{i-1} + (k-c_i-b_i) u_i + b_i u_{i+1} = \theta u_i$,
where $c_i = i^2$ and $b_i = (m-1-i)(n-i)$ and $k = (m-1)n$.
Define $v_i = iu_{i-1} + (m-i)u_i$ ($1 \le i \le n$).
Then $c'_i v_{i-1} + (k-c'_i-b'_i) v_i + b'_i v_{i+1} = \theta v_i$,
where $c'_i = i(i-1)$ and $b'_i = (n-i)(m-i)$.
It follows that $Ev = \theta v$. We have $v \ne 0$ for $\theta \ne -n$.
\qed

\noindent
{\bf Remark.}\quad The Johnson graph $J(m,n)$ has eigenvalues
$(n-i)(m+1-n-i)-n$ with multiplicity ${m \choose i} - {m \choose i-1}$
$(0 \le i \le \min(n,m-n))$. It is not wrong to say that it has these
eigenvalues and multiplicities for $0 \le i \le n$ since by convention
the multiplicities of an eigenvalue are added, and eigenvalues with
multiplicity 0 are no eigenvalues. For example, $J(5,4)$ has
spectrum $4^1$ $(-1)^4$, but it also has spectrum
$4^1$ $(-1)^4$ $(-4)^5$ $(-5)^0$ $(-4)^{-5}$, where multiplicities
are written as exponents.
The Johnson graph $J(m+n-1,n)$ has eigenvalues $(m-i)(n-i)-n$
with multiplicity ${m+n-1 \choose i} - {m+n-1 \choose i-1}$
$(0 \le i \le \min(n,m-1))$.

\section{The common part of the spectra of ${\rm SR}(m,n)$ and $J(m+n-1,n)$}
Both ${\rm SR}(m,n)$ and $J(m+n-1,n)$ have  ${m+n-1 \choose n}$ vertices.
Both have valency $n(m-1)$.
These graphs resemble each other and have a large part of their spectrum
in common.
Let $m,n > 0$.

\begin{Proposition} \label{commonev}
The graphs ${\rm SR}(m,n)$ and $J(m+n-1,n)$ have equitable partitions
with the same quotient matrix $E$, where $E$ has eigenvalues
$(n-i)(m-i)-n$ with multiplicity ${m \choose i}$ for $0 \le i \le n-1$,
and multiplicity ${m \choose n} - 1$ for $i = n$.
In particular, the spectrum of $E$ is a common part of the spectrum of
${\rm SR}(m,n)$ and that of $J(m+n-1,n)$.
\end{Proposition}

\noindent
That is, $E$ has eigenvalue $(n-i)(m-i)-n$ with multiplicity
${m \choose i}$ for $0 \le i \le \min(m,n)-1$,
and if $n < m$ also eigenvalue $-n$ with multiplicity ${m \choose n} - 1$.

\medskip\noindent
\Proof
Partition the vertex set of ${\rm SR}(m,n)$ into $\sum_{i=1}^n {m \choose i}$
parts, where each part consists of the vertices with fixed support $S$
(of weight $i$).
Partition the vertex set of $J(m+n-1,n)$ into $\sum_{i=1}^n {m \choose i}$
parts, where each part consists of the vertices with fixed support $S$
(of weight $i$) in the first $m$ coordinates.
The size of part $S$, where $|S|=i$, is ${n-1 \choose n-i}$ in both cases.

It is straightforward to determine the numbers $e_{ST}$ for both graphs.
For ${\rm SR}(m,n)$, note that if we restrict attention
to vertices with support $S$, where $|S|=i$, then our graph
becomes a copy of ${\rm SR}(i,n-i)$, as one sees by subtracting 1
from all entries. The result is
$$
e_{ST} = \left\{ \begin{tabular}{ll}
$(i-1)(n-i)$ & if $|S|=i$, $S = T$, \\
$i-1$ & if $|S|=i$, $|T|=i-1$, $S \supset T$, \\
$n-i$ & if $|S|=i$, $|T|=i+1$, $S \subset T$, \\
1 & if $|S|=|T|=i$ and $S,T$ differ in two places, \\
0 & otherwise.
\end{tabular} \right.
$$
in both cases.
It follows that our partitions are equitable with the same
quotient matrix $E$. We may conclude that $J(m+n-1,n)$ and
${\rm SR}(m,n)$ have the $\sum_{i=1}^n {m \choose i}$ eigenvalues
of the matrix $E$ in common.

\medskip
Claim: these eigenvalues are
$(n-i)(m-i)-n$ with multiplicity ${m \choose i}$ for $0 \le i < n$,
and multiplicity ${m \choose n} - 1$ for $i = n$.
These are the eigenvalues of $J(m+n-1,n)$, so we need only
confirm the multiplicities.

\medskip
Let $W_{ij}$ be the (symmetrized) inclusion matrix of $i$-subsets
against $j$-subsets in a $v$-set. Then $W_{ij} = W_{ji}^\top$.
For $h \le i \le j$ we have $W_{hi}W_{ij} = {j-h \choose i-h} W_{hj}$.
Also $W_{ii} = I$ and
$W_{i,j+1}W_{j+1,j} = W_{i,i-1}W_{i-1,j} + (v-i-j)W_{i,j}$
for $j \ge i$. (Note that $W_{i,-1}$ has no columns, and
$W_{-1,j}$ has no rows, so that $W_{i,-1}W_{-1,j} = 0$.)

For $0 \le i \le n$, the $i$-th eigenspace of $J(v,n)$ is spanned by
vectors $W_{ni} w$ where $w$ is indexed by the $i$-subsets of the $v$-set
and $W_{i-1,i} w = 0$ if $i > 0$.
(Indeed, the adjacency matrix of $J(v,n)$ is $A = W_{n,n-1} W_{n-1,n} - nI$,
and $W_{n,n-1} W_{n-1,n} W_{ni} w =
W_{n,n-1} (W_{n-1,i-1} W_{i-1,i} + (v-n-i+1) W_{n-1,i}) w =
(v-n-i+1)(n-i) W_{ni} w$,
so that $J(v,n)$ has eigenvalue $(v-n-i+1)(n-i)-n$ with
multiplicity ${v \choose i} - {v \choose i-1}$.)

\medskip
Let $v = m+n-1$, and consider the vectors $W_{ni} w$ that are invariant
under ${\rm Sym}(n-1)$ acting on the final $n-1$ coordinates.
They form a space of dimension $\sum_{j=0}^i {m \choose j}$
if $i < n$, and $\sum_{j=1}^n {m \choose j}$ if $i = n$.
The subspace of such vectors satisfying $W_{i-1,i} w = 0$
has codimension $\sum_{j=0}^{i-1} {m \choose j}$:
each $(i-1)$-set $T$ imposes a restriction $(W_{i-1,i} w)_T = 0$,
and such restrictions are equivalent when the $T$'s are in the
same ${\rm Sym}(n-1)$-orbit.
It follows that the corresponding eigenspace of $E$ has
dimension ${m \choose i}$ when $i < n$. For $i=n$ this same
computation becomes $(\sum_{j=1}^n {m \choose j}) -
(\sum_{j=0}^{n-1} {m \choose j}) = {m \choose n} - 1$.
\qed

\section{Automorphism group}
Clearly, $G = {\rm Sym}(m)$ acts as a group of automorphisms
on $\Gamma$, permuting the coordinate positions.
This action is faithful, unless $n=0$, $m > 1$ when the full group has order 1.
When $m=2$, the graph is $K_{n+1}$ and the full group is ${\rm Sym}(n+1)$.
When $n=2$, the graph is $T(m+1) \isom J(m+1,2)$ with full group
${\rm Sym}(m+1)$ for $m \ne 1,3$ and $2^3.{\rm Sym}(3)$ for $m=3$.

\begin{Proposition} \label{automgp}
Let $m,n > 2$.
If $n = 3$, then ${\rm Aut} (\Gamma) \isom {\rm Sym}(m).2$,
where the additional factor $2$ interchanges the digits $1$ and $2$
in each vector with a coordinate $2$.
If $n > 3$, then ${\rm Aut} (\Gamma) \isom {\rm Sym}(m)$.
\end{Proposition}
\Proof
Classify the vertices $x$ according to the number $\lambda_{xy}$
of common neighbors of $x$ and $y$, for all neighbors $y$ of $x$.

We always have $\lambda_{xy} \le m+n-3$.
(Indeed, look at the common neighbors of $x = (a,b,c,...)$ and
$y = (a+d,b-d,c,...)$. We find $(a',b',c,...)$ with $a'+b' = a+b$
($a+b-1$ choices), and $(a,b-d,c+d,...)$ ($m-2$ choices),
and $(a+d,b,c-d,...)$ where the number of choices is the number of $c$
not less than $d$. If $a+b=n$, there are no such $c$. Otherwise this
number is maximal when $d = 1$ and all nonzero $c$ equal 1,
and then equals $n-a-b$. So $\lambda_{xy} \le a+b-1+m-2+n-a-b = m+n-3$.)

Note for later use the structure of the graph $\Lambda(x,y)$ induced by
the common neighbors of $x$ and $y$. It is $K_{a+b-1} + K_{m-2} + K_g$,
where $g$ is the number of $c$ (common coordinates of $x$ and $y$)
not less than $d$.

If $\lambda_{xy} = m+n-3$ for all $n(m-1)$ neighbors $y$ of $x$
(and $m > 2$) then $x$ either has a unique nonzero coordinate $n$,
or only has coordinates 0, 1. Thus, we can recognize this set of
$m + {m \choose n}$ vertices. The induced subgraph (for $n > 2$)
is isomorphic to $K_m + J(m,n)$.

We see that $\Gamma$ determines $m+n-3$ and also the pair
$\{m, {m \choose n}\}$. Now $m$ is the smallest element element
of the pair distinct from~0,~1, so we find $m$ and $n$.

Suppose first that $n \ne m-1$. Then we recognized the set $S$ of
vectors with a unique nonzero coordinate. At distance $i$ from $S$
lie the vectors with precisely $i+1$ nonzero coordinates, and the
positions of the nonzero coordinates of a vector $u$ are determined
by the set of nearest vertices in $S$. We show by induction on $m$
that all vertex labels are determined. If a vertex $(a,b,...)$
has at least two nonzero coordinates, and $m > 3$, then its neighbor
$(0,a+b,...)$ lies in the ${\rm SR}(m-1,n)$ on the vertices with first
coordinate zero, and by induction $a+b$ is determined.
If it has at least three nonzero coordinates: $(a,b,c,...)$,
then each of $a+b$, $a+c$, $b+c$ is determined, and hence also $a,b,c$.
If it has precisely two nonzero coordinates: $(a,n-a,0,...)$ then it
has neighbors $(a-i,n-a,i,0,...)$ ($1 \le i \le a-1$) and
$(a,n-a-j,j,0,...)$ ($1 \le j \le n-a-1$) of which all coordinates
are known, and $a$ and $n-a$ follow unless $\{a,n-a\} = \{1,2\}$.
This settles all claims when $m > 3$, $n \ne m-1$.

If $m = 3$, $n \ge 3$, we recall that the common neighbors of
vertices $x$ and $y$ induce $K_{a+b-1} + K_{m-2} + K_g$,
where $m-2=1$ and $g \le 1$, so that $a+b$ can be recognized directly
when $a+b \ge 3$. But we also know the zero pattern, so can also
recognize $a,b$ when $a+b \le 2$. This determines all for $m = 3$.

Finally, if $m = n+1 \ge 4$, we have to distinguish the copy of $K_m$
on the vectors of shape $n^1 0^n$ from that on the vectors of shape $1^n 0$.
Both sets have $m{m-1 \choose 2}$ neighbors, but if
$x=(n,0,...)$ and $y=(n-a,a,...)$ then $\Lambda(x,y) \isom 2K_{n-1}$,
while $\Lambda(x,y) \isom K_{n-1}+K_{n-2}+K_1$ if $x=(1,...,1,1,0)$
and $y=(2,1,...,1,0,0)$.
This settles all cases.
\qed

\section{Diameter}
\begin{Proposition} \label{diameter}
Let $m > 0$. Then the diameter of ${\rm SR}(m,n)$ is $\min (m-1,n)$.
\end{Proposition}
\Proof
The diameter is at most $m-1$, since one can walk from one
vertex to another and decrease the number of different coordinates
by at least one at each step.
The diameter is also at most $n$, since one can walk from one vertex
to another and decrease the sum of the absolute values of the coordinate
differences by at least two at each step.
If $m > n$, then $(0,...,0,n)$ and $(1,...,1,0,...,0)$
show that the diameter is at least $n$.
If $m \le n$, then $(0,...,0,n)$ and $(1,...,1,n-m+1)$
show that the diameter is at least $m-1$.
\qed

\section{Maximal cliques and local graphs}
We classify the cliques (complete subgraphs), and find the maximal ones.
We also examine the structure of the local graphs of $\Gamma$.

\begin{Lemma}
Cliques $C$ in ${\rm SR}(m,n)$ are of three types:

1. All adjacencies are $(j,k)$-adjacencies for fixed $j,k$.
Now $|C| \le n+1$.

2. $C = \{ x{+}ae_i \mid i \in I \}$, where $a \in \Nn$, $1 \le a \le n$,
$x \in \Nn^m$ with $\sum x_i = n-a$, and $I \subseteq \{1,\ldots,m\}$.
Now $|C| \le m$.

3. $C = \{x{-}ae_i \mid i \in I \}$, where $a \in \Nn$, $a \ge 1$,
$x \in \Nn^m$ with $\sum x_i = n+a$, $I \subseteq \{1,\ldots,m\}$,
and $x_i \ge a$ for $i \in I$.
Now $|C| \le m$.
\end{Lemma}
\Proof
Suppose $u,v,w$ are pairwise adjacent, not all $(j,k)$-adjacent
for the same pair $(j,k)$. Then $u,v$ are $(i,j)$-adjacent,
$u,w$ are $(i,k)$-adjacent, and $v,w$ are $(j,k)$-adjacent,
for certain $i,j,k$.
Now $u_k=v_k$, $u_j=w_j$ and $v_i=w_i$, so that
$u = x+ae_i$, $v = x+ae_j$, $w = x+ae_k$, where $a > 0$ or $a < 0$.
\qed

\begin{Proposition} \label{cliquesz}
Let $m > 1$ and $n > 0$.
Then the size of the largest clique in ${\rm SR}(m,n)$ is $\max(m,n+1)$.
\end{Proposition}
\Proof
For $n > 0$ the $m$ vectors $ne_i$ are distinct and
mutually adjacent, forming an $m$-clique.
And for $m > 1$ the $n+1$ vectors $ae_1+(n-a)e_2$ $(0 \le a \le n)$
form an $(n+1)$-clique. Conversely, no larger cliques occur, as we
just saw.
\qed

\medskip\noindent
Fix a vertex $u$ of ${\rm SR}(m,n)$. We describe the structure
of the local graph of $u$, that is the graph induced by ${\rm SR}(m,n)$
on the set $U$ of neighbors of $u$. If $vw$ is an edge in this
local graph, then $uvw$ is a clique in ${\rm SR}(m,n)$, so we can
invoke the above classification.

\begin{Lemma}
(i) Any two adjacent vertices $u,v$ uniquely determine three cliques
$C_1,C_2,C_3$ where $C_i$ is of type $i$ and
$C_i \cap C_j = \{u,v\}$ for distinct $i,j \in \{1,2,3\}$,
and $C_1 \cup C_2 \cup C_3$ contains all common neighbors of $u$ and $v$.

(ii) Fix a vertex $u$. For each $i \in \{1,2,3\}$ the
cliques on $u$ of type $i$ form (after removal of $u$)
a partition of the set $U$ of neighbors of $u$.
Each edge in $U$ is contained in a unique such clique,
and hence has a unique type.

(iii) ${\rm SR}(m,n)$ does not contain an induced $K_{1,1,4}$.
\end{Lemma}
\Proof
Let $C_i = C_i(u,v)$ $(i=1,2,3)$ be the unique largest clique on $\{u,v\}$
of type $i$ $(i=1,2,3)$. \qed

The set $U$ has a partition into ${m \choose 2}$ cliques of type 1,
where the $(j,k)$-clique has size $u_j+u_k$.
(We check that $\sum_{j,k} u_j+u_k = (m-1) \sum_j u_j = (m-1)n$.)\\
The set $U$ has a partition into $n$ cliques of type 2, each of size $m-1$.
Finally, $U$ has a partition into cliques of type 3.
(If $v = u-ae_j+ae_k$ is a neighbor of $u$, then
$C_3(u,v) = \{ x-ae_i \mid i \in I \}$, where $x = u+ae_k$
and $I = \{ i \mid 1 \le i \le m,~ x_i \ge a \}$.)

\begin{Lemma} \label{noK113}
Let $m,n \ge 3$, and fix a vertex $u$.
Each neighbor $v$ of $u$ is contained in at most two maximal cliques
precisely when $u$ has only one nonzero coordinate.
\end{Lemma}
\Proof
Suppose each point $v$ of $U$ is covered by at most two maximal cliques.
Then one of the cliques of types 1 or 3 on $v$ in $U$ has size 1.
This means that whenever $u_j+u_k \ge 2$, we have $u_i=0$ for $i \ne j,k$.
If $u_j \ge 2$ this means that $u$ has only one nonzero coordinate.
If $u_j = u_k = 1$ this means that $n = 2$. \qed

Suppose $m,n \ge 3$. We see that we can retrieve $V_1$
as the set of vertices that are locally locally the union of
two cliques.

\section{Independence number}

It is known that ${\rm SR}(3,n)$ has independence number
$\alpha(3,n) = \lfloor (2n+3)/3 \rfloor$
(see \cite[Problem 252]{VGL}, \cite{NL}, \cite{BPS}).
We determine $\alpha(m,3)$.

\begin{Proposition}
$$
\alpha(m,3) = \left\{ \begin{array}{ll}
\frac16 (m+1)(m+2) & {\rm for}~ m \equiv \pm 1 ~ {\rm mod} ~ 6, \\
\frac16 m(m+3) & {\rm for}~ m \equiv 3 ~ {\rm mod} ~ 6, \\
\frac16 m(m+2) & {\rm for}~ m \equiv 0,4 ~ {\rm mod} ~ 6, \\
\frac16 (m^2+2m-2) &  {\rm for}~ m \equiv 2 ~ {\rm mod} ~ 6. \\
\end{array}\right.
$$
\end{Proposition}
\Proof
A trivial upper bound is
$$\alpha(m,3) \le
m + \lfloor \frac13 ( m \lfloor \frac12 (m-3) \rfloor + 1) \rfloor.$$

(As follows: count edges in $K_m$ covered by vertices of the independent set.
No edge is covered twice because the corresponding vertices would be adjacent.
A vertex $3e_i$ (singleton) covers no edge.
A vertex $2e_i+e_j$ (pair) covers the edge $ij$.
A vertex $e_i+e_j+e_k$ (triple) covers the edges $ij$ and $ik$ and $jk$.
Nonadjacent vertices $3e_h$ and $2e_i+e_j$ have $h \ne i,j$.
Nonadjacent vertices $2e_i+e_j$ and $2e_k+e_l$ have
$i \ne k$ and $j \ne l$ and $\{i,j\} \ne \{k,l\}$.
It follows that there are at most $m$ vertices of the forms
$3e_h$ and $2e_i+e_j$. Since a triple takes 3 edges, and a pair only 1,
and a singleton 0, we find an upper bound by assuming that there are
1 singleton and $m-1$ pairs (2 on each non-singleton). That leaves
$m-3$ edges on each non-singleton, and 2 more on the singleton, for
a maximum of $\frac12 (m-3)$ triples on each non-singleton, and
1 more on the singleton. Since each triple is counted thrice,
there are at most $\frac13 ( m \lfloor \frac12 (m-3) \rfloor + 1)$
triples.)

Separating the cases for $m$ (mod 6) yields precisely the values claimed,
so they are upper bounds. But examples reaching the bounds can be
constructed from Steiner triple systems.

If $m \equiv \pm 1$~(mod~6) then $\alpha(m,3) \ge \frac16 (m+2)(m+1)$.
(Take an STS($m+2$), and delete two points $x,y$. The graph of
noncovered edges has valency 2, so is a union of cycles.
Direct each cycle. Pick $2e_i+e_j$ for each directed edge $(i,j)$.
Let $xyz$ be a block. Pick $3e_z$. This is a coclique
of the indicated size that reaches the trivial upper bound.)

If $m \equiv 3$ (mod 6), then $\alpha(m,3) \ge \frac16 m(m+3)$.
(Take an STS($m$) with parallel class (e.g., a KTS($m$))
and view the triples in the parallel class as directed cycles.)

If $m \equiv 0,2,4$ (mod 6), then the claimed value is obtained
by shortening an example for $m+1$.
\qed

\section{Cospectral mates}
For $m \le 2$ or $n \le 2$, the graph ${\rm SR}(m,n)$ is
complete or triangular, and hence determined by its spectrum,
except in the case of $m=7$, $n=2$ where it is isomorphic to
the triangular graph $T(8)$, and cospectral with the three Chang graphs
(cf.~\cite{Chang59,Chang60}).
The graph ${\rm SR}(3,3)$ is 6-regular on 10 vertices, and we find
that its complement is cubic with spectrum $3^1$ $2^1$ $1^3$ $(-1)^2$ $(-2)^3$.
All integral cubic graphs are known, and ${\rm SR}(3,3)$ is
uniquely determined by its spectrum, cf.~\cite{BH}, \S3.8.
We give some further cases where ${\rm SR}(m,n)$ is not determined
by its spectrum.

\begin{Proposition} \label{cospectral}
The graph ${\rm SR}(m,n)$ is not determined by its spectrum
when (a)~$m=4$ and $n \ge 3$, or (b)~$n = 3$ and $m \ge 4$.
\end{Proposition}
\Proof
Apply Godsil-McKay switching (cf.~\cite{GM}; \cite{BH}, 1.8.3, 14.2.3).
Switch with respect to a 4-clique $B$ such that every vertex outside $B$
is adjacent to 0, 2 or 4 vertices inside.
If $m=4$, take $B = \{ n000, 0n00, 00n0, 000n \}$.
If $n=3$, $m \ge 2$, take $B = \{ ae_1+be_2 \mid a+b=3 \}$.
In both cases every vertex outside $B$ is adjacent to 0 or 2 vertices inside.
The switching operation preserves all edges and nonedges, except
that it changes adjacency for pairs $bc$ with $b \in B$, $c \notin B$,
and $c$ adjacent to 2 vertices of $B$,
turning edges (resp.~nonedges) into nonedges (resp.~edges).
The resulting graph has the same spectrum. We show that it is
nonisomorphic to ${\rm SR}(m,n)$ for $m=4$, $n \ge 3$
and for $n = 3$, $m \ge 4$.

In the former case, $B = V_1$. If switching does not change the isomorphism
type, then $B$ must remain the $V_1$ of the new graph (since it is
a single orbit of size $m$ contained in $V_1 \cup V_2$).
But after switching the common neighbours of
$n000$ and $0n'10$ (with $n'=n-1$) include the pairwise
nonadjacent $0n'01$, $01n'0$, $001n'$, contradicting Lemma \ref{noK113}.

In the latter case, $B = \{ 3000..,~ 2100..,~ 1200..,~ 0300.. \}$.
After switching, $3000..$ and $0300..$ are still in $V_1$
since their local graphs are $3 \times (m-1)$ grids.
But $2100..$ and $1200..$ are not, since the common neighbors of
$2100..$ and $0030..$ include the pairwise nonadjacent
$1200.., 1020.., 0210..$.
And $1110..$ is not, since the common neighbors of
$1110..$ and $0120..$ include the pairwise nonadjacent
$3000..,1020..,0111..$.
So, there is no candidate for $V_1$.
\qed

There are at least 336 pairwise nonisomorphic graphs with spectrum
$9^1$ $3^4$ $1^3$ $(-1)^6$ $(-3)^6$, namely
${\rm SR}(4,3)$ and the three graphs obtained by Godsil-McKay switching
with respect to the 4-cliques $\{ 3000, 0300, 0030, 0003 \}$,
$\{ 0111, 1011,\maysplit 1101,\maysplit 1110 \}$ and
$\{ 3000, 2100, 1200, 0300 \}$,
and 332 further graphs obtained by repeated switching
w.r.t. regular subgraphs of size 4.

\section{The eigenspace of the smallest eigenvalue}

Fix $\pi \in {\rm Sym}(m)$, and let
$a_i = \#\{ j \mid i < j ~{\rm and}~\pi_i > \pi_j \}$
for $1 \le i \le m$.
Then $a = (a_i)$ is a vertex of ${\rm SR}(m,n)$ when $n$ is the number
of inversions of $\pi$.

Say that $\sigma \in {\rm Sym}(m)$ is
$\pi$-{\em admissible} if $a_i + i - \sigma_i \ge 0$ for $1 \le i \le m$.
Let ${\rm Adm}(\pi)$ be the set of $\pi$-admissible permutations,
and define $x(\sigma)$ by $x(\sigma)_i = a_i+i-\sigma_i$.
Then $\sigma \in {\rm Adm}(\pi)$ if and only if
$x(\sigma)$ is a vertex of ${\rm SR}(m,n)$.

\begin{Theorem} {\rm (Martin \& Wagner \cite{MW}, Thm.~3.8)} \label{parp}
For each $\pi \in {\rm Sym}(m)$ with $n$ inversions, let
$$F_\pi = \sum_{\sigma \in {\rm Adm}(\pi)} {\rm sgn}(\sigma) e_{x(\sigma)}.$$
Then each $F_\pi$ is an eigenvector of ${\rm SR}(m,n)$ with eigenvalue~$-n$,
and the $F_\pi$ are linearly independent.
\end{Theorem}

\begin{Theorem} {\rm (Martin \& Wagner \cite{MW}, Prop.~3.1)}
For $p,w \in \Rr^m$ such that $p+\sigma(w)$ for $\sigma \in {\rm Sym}(m)$
are pairwise distinct vertices of ${\rm SR}(m,n)$, define
$$F_{p,w} = \sum_{\sigma \in {\rm Sym}(m)} {\rm sgn}(\sigma)
e_{p+\sigma(w)}.$$
Then each $F_{p,w}$ is an eigenvector of ${\rm SR}(m,n)$ with eigenvalue
$-{m \choose 2}$, and for fixed $w$, the collection of all such $F_{p,w}$
is linearly independent.
\end{Theorem}

Picking $w = \frac12(1-m,3-m,\ldots,m-3,m-1)$ yields the lower bound
already mentioned earlier: the multiplicity of the eigenvalue
$-{m \choose 2}$ is at least ${n-{m-1 \choose 2} \choose m-1}$.

For the eigenvalue $-n$ it follows that its multiplicity is
at least the number of elements in ${\rm Sym}(m)$ with precisely
$n$ inversions, and one conjectures that equality holds.

\medskip
The proof of Theorem \ref{parp} shows that for each $\pi \in {\rm Sym}(m)$
with $n$ inversions, the set
$X_\pi = \{ x(\sigma) \mid \sigma \in {\rm Adm}(\pi) \}$
induces a bipartite subgraph of ${\rm SR}(m,n)$ that is regular
of valency $n$.

\medskip
One may wonder what graphs $\Gamma(m,n,\pi)$ occur as induced subgraph
on such a subset $X_\pi$ of the vertex set of ${\rm SR}(m,n)$.
Given $\Gamma(m,n,\pi)$, one can find a $\pi'$ such that
$\Gamma(m+1,n,\pi') \isom \Gamma(m,n,\pi)$, and a $\pi''$ such that
$\Gamma(m+2,n+1,\pi'') \isom \Gamma(m,n,\pi) \times K_2$, where
$\times$ denotes Cartesian product.

\medskip
If $n = {m \choose 2}$, then there is a unique permutation
$\pi_0 \in {\rm Sym}(m)$ with $n$ inversions. It has
$a = (m-1,m-2,\ldots,0)$, and all $n!$ permutations are $\pi_0$-admissible.
The resulting graph $\Gamma(m,{m \choose 2},\pi_0)$ has as vertices all
permutations of $(m-1,m-2,\ldots,0)$, where two permutations are adjacent
when they differ by a transposition. In other words, this is the Cayley
graph ${\rm Cay}({\rm Sym}(m),T)$, where $T$ is the set of transpositions
in ${\rm Sym}(m)$.

\begin{Proposition}
For any $\Gamma(m,n,\pi)$, where $m > 2n$, there is an isomorphic
$\Gamma(2n,n,\pi')$.
\end{Proposition}

It follows that classifying all $\Gamma(m,n,\pi)$
for fixed $n$ is a finite job. Let $Q_k$ denote the $k$-cube.
Using Sage we find for $n=1$ that only $Q_1$ occurs,
for $n=2$ that only $Q_2$ occurs,
for $n=3$ that only $K_{3,3}$ and $Q_3$ occur,
and for $n=4$ that only $K_{3,3} \times K_2$ and $Q_4$ occur.
For larger $n$ one finds more complicated shapes.

It was conjectured in \cite{MW} that all graphs $\Gamma(m,n,\pi)$
have integral spectrum.

\section{Spectra for small $m$ or $n$}
If we fix a small value of $n$, we find a nice spectrum
(eigenvalues and multiplicities are polynomials in $m,n$).
If we fix a small value of $m \ge 3$, we get a messy result
(also congruence conditions play a r\^ole).
Below, multiplicities are written as exponents.

\medskip
For $n=0$, the spectrum is $0^1$.

\medskip
For $n=1$, the spectrum is $(m-1)^1$, $(-1)^{m-1}$.

\medskip
For $n=2$, the spectrum is $(2m-2)^1$, $(m-3)^m$, $(-2)^{\rm rest}$. 

\medskip
For $n=3$, the spectrum is
$(3m-3)^1$, $(2m-5)^m$, $(m-3)^{m-1}$, $(m-5)^{m \choose 2}$,
$(-3)^{\rm rest}$. (See below.)

\medskip
For $n=4$, the spectrum is
$(4m-4)^1$, $(3m-7)^m$, $(2m-5)^m$, $(2m-8)^{m \choose 2}$,\\
$(m-4)^{{m \choose 2}-1}$, $(m-6)^{m \choose 2}$,
$(m-7)^{m \choose 3}$, $(-4)^{\rm rest}$. (See below.)

\medskip
For $n=5$, the spectrum may be
$(5m-5)^1$, $(4m-9)^m$, $(3m-7)^m$,\\
$(3m-11)^{m \choose 2}$, $(2m-5)^{m-1}$,
$(2m-7)^{m \choose 2}$, $(2m-9)^{m \choose 2}$,
$(2m-11)^{m \choose 3}$, $(m-5)^{{m \choose 3}-1}$,\\
$(m-6)^{m(m-2)}$,
$(m-8)^{2{m \choose 3}}$, $(m-9)^{m \choose 4}$, $(-5)^{\rm rest}$.

\bigskip
For $m=1$, the spectrum is $0^1$.

\medskip
For $m=2$, the spectrum is $n^1$, $(-1)^n$.

\medskip
For $m=3$, the spectrum is
$(2n)^1$,
$b^3$ (for all $b \in \Zz$ with $-2 \le b \le n-2$),
$(-3)^{n-1 \choose 2}$,
except for $(a-1)^3$ and $a^1$ (if $n=2a+3$)
or $a^3$ and $(a-1)^1$ (if $n=2a+4$) (Martin \& Wagner \cite{MW}).

\medskip
For $m=4$, we give some values in Table \ref{table4}.

Let $a^m \downarrow b$ denote sequence of eigenvalues and multiplicities
found as follows: the eigenvalues are the integers $c$ with
$a \ge c \ge b$, where the first multiplicity is $m$,
and each following multiplicity is 2 larger for even $c$,
and 10 larger for odd $c$.
Now the conjectured spectrum of ${\rm SR}(4,n)$, $n \ge 6$, $n \ne 7$
consists of

(i) $(3n)^1$,

\smallskip
(ii) $b^4$ for all odd integers $b$, where $2n-3 \ge b \ge n-1$,

\smallskip
(iii)

\newcommand\Tstrut{\rule{0pt}{10pt}}
\noindent
\begin{tabular}{@{\,}l|l@{}}
\hline
\Tstrut$n = 2r$ & $(n-4)^{3n-1}$, $(n-6)^6$,
 $(n-7)^{16} \downarrow (n-8)/2$ \\
$n = 2r+1$ & $(n-2)^3$, $(n-4)^{3n-3}$, $(n-6)^9$,
 $(n-7)^{12} \downarrow (n-7)/2$ \\
\hline
\end{tabular}

\smallskip
(iv) for $q = \lceil n/3-4 \rceil$:

\noindent
\begin{tabular}{@{\,}l|l@{}}
\hline
\Tstrut$n = 4s$ & $(2s-5)^{3n-12}$, $(2s-6)^{3n-26} \downarrow q$ \\
$n = 4s+1$ & $(2s-4)^{3n-7}$, $(2s-5)^{3n-21}$,
 $(2s-6)^{3n-23} \downarrow q$ \\
$n = 4s+2$ & $(2s-4)^{3n-16}$, $(2s-5)^{3n-22} \downarrow q$ \\
$n = 4s+3$ & $(2s-3)^{3n-3}$, $(2s-4)^{3n-25}$,
 $(2s-5)^{3n-19} \downarrow q$ \\
\hline
\end{tabular}

\smallskip
(v) if $n \equiv 0$ (mod 3) one additional eigenvalue $n/3-4$,

\smallskip
(vi)

\noindent
\begin{tabular}{@{\,}l|l@{}}
\hline
\Tstrut$n = 6t$ & $(2t-5)^{4n-12}$, $(2t-6)^{4n-16}$,
 $(2t-7)^{4n-16} \downarrow (-5)$ \\
$n = 6t+1$ & $(2t-4)^{4n-32}$, $(2t-5)^{4n-7}$,
 $(2t-6)^{4n-20}$, $(2t-7)^{4n-14} \downarrow (-5)$ \\
$n = 6t+2$ & $(2t-4)^{4n-24}$, $(2t-5)^{4n-8}$,
 $(2t-6)^{4n-21}$, $(2t-7)^{4n-12} \downarrow (-5)$ \\
$n = 6t+3$ & $(2t-4)^{4n-16}$, $(2t-5)^{4n-12}$,
 $(2t-6)^{4n-20} \downarrow (-5)$ \\
$n = 6t+4$ & $(2t-3)^{4n-28}$, $(2t-4)^{4n-11}$,
 $(2t-5)^{4n-16}$, $(2t-6)^{4n-18} \downarrow (-5)$ \\
$n = 6t+5$ & $(2t-3)^{4n-20}$, $(2t-4)^{4n-12}$,
 $(2t-5)^{4n-17}$, $(2t-6)^{4n-16} \downarrow (-5)$ \\
\hline
\end{tabular}

\smallskip
(vii) finally $(-6)^{n-3 \choose 3}$.

\medskip
For example, $-5$ has multiplicity $6n-28$.

\begin{table}[ht]
\begin{center}
\begin{tabular}{| l | p{.9\textwidth} |}
\hline
$n$ & Spectrum of $SR(4,n)$ \\ \hline
0 & $0^1$ \\
1 & $3^{1}$ $(-1)^{3}$ \\
2 & $6^{1}$ $1^{4}$ $(-2)^{5}$ \\
3 & $9^{1}$ $3^{4}$ $1^{3}$ $(-1)^{6}$ $(-3)^{6}$ \\
4 & $12^{1}$ $5^{4}$ $3^{4}$ $0^{11}$ $(-2)^{6}$ $(-3)^{4}$ $(-4)^{5}$ \\
5 & $15^{1}$ $7^{4}$ $5^{4}$ $3^{3}$ $1^{12}$ $(-1)^{9}$ $(-2)^{8}$ $(-3)^{4}$ $(-4)^{8}$ $(-5)^{3}$ \\
6 & $18^{1}$ $9^{4}$ $7^{4}$ $5^{4}$ $2^{17}$ $0^{6}$ $(-1)^{16}$ $(-2)^{3}$ $(-3)^{12}$ $(-4)^{8}$ $(-5)^{8}$ $(-6)^{1}$ \\
7 & $21^{1}$ $11^{4}$ $9^{4}$ $7^{4}$ $5^{3}$ $3^{18}$ $1^{9}$ $0^{12}$ $(-1)^{18}$ $(-3)^{21}$ $(-4)^{8}$ $(-5)^{14}$ $(-6)^{4}$ \\
8 & $24^{1}$ $13^{4}$ $11^{4}$ $9^{4}$ $7^{4}$ $4^{23}$ $2^{6}$ $1^{16}$ $0^{18}$ $(-1)^{12}$ $(-2)^{8}$ $(-3)^{24}$ $(-4)^{11}$ $(-5)^{20}$ $(-6)^{10}$ \\
9 & $27^{1}$ $15^{4}$ $13^{4}$ $11^{4}$ $9^{4}$ $7^{3}$ $5^{24}$ $3^{9}$ $2^{12}$ $1^{22}$ $0^{20}$ $(-1)^{7}$ $(-2)^{20}$ $(-3)^{24}$ $(-4)^{16}$ $(-5)^{26}$ $(-6)^{20}$ \\
10 & $30^{1}$ $17^{4}$ $15^{4}$ $13^{4}$ $11^{4}$ $9^{4}$ $6^{29}$ $4^{6}$ $3^{16}$ $2^{18}$ $1^{28}$ $0^{14}$ $(-1)^{12}$ $(-2)^{29}$ $(-3)^{24}$ $(-4)^{22}$ $(-5)^{32}$ $(-6)^{35}$ \\
11 & $33^{1}$ $19^{4}$ $17^{4}$ $15^{4}$ $13^{4}$ $11^{4}$ $9^{3}$ $7^{30}$ $5^{9}$ $4^{12}$ $3^{22}$ $2^{24}$ $1^{30}$ $0^{8}$ $(-1)^{24}$ $(-2)^{32}$ $(-3)^{27}$ $(-4)^{28}$ $(-5)^{38}$ $(-6)^{56}$ \\
12 & $36^{1}$ $21^{4}$ $19^{4}$ $17^{4}$ $15^{4}$ $13^{4}$ $11^{4}$ $8^{35}$ $6^{6}$ $5^{16}$ $4^{18}$ $3^{28}$ $2^{30}$ $1^{24}$ $0^{11}$ $(-1)^{36}$ $(-2)^{32}$ $(-3)^{32}$ $(-4)^{34}$ $(-5)^{44}$ $(-6)^{84}$ \\
13 & $39^{1}$ $23^{4}$ $21^{4}$ $19^{4}$ $17^{4}$ $15^{4}$ $13^{4}$ $11^{3}$ $9^{36}$ $7^{9}$ $6^{12}$ $5^{22}$ $4^{24}$ $3^{34}$ $2^{32}$ $1^{18}$ $0^{20}$ $(-1)^{45}$ $(-2)^{32}$ $(-3)^{38}$ $(-4)^{40}$ $(-5)^{50}$ $(-6)^{120}$ \\
14 & $42^{1}$ $25^{4}$ $23^{4}$ $21^{4}$ $19^{4}$ $17^{4}$ $15^{4}$ $13^{4}$ $10^{41}$ $8^{6}$ $7^{16}$ $6^{18}$ $5^{28}$ $4^{30}$ $3^{40}$ $2^{26}$ $1^{20}$ $0^{32}$ $(-1)^{48}$ $(-2)^{35}$ $(-3)^{44}$ $(-4)^{46}$ $(-5)^{56}$ $(-6)^{165}$ \\
15 & $45^{1}$ $27^{4}$ $25^{4}$ $23^{4}$ $21^{4}$ $19^{4}$ $17^{4}$ $15^{4}$ $13^{3}$ $11^{42}$ $9^{9}$ $8^{12}$ $7^{22}$ $6^{24}$ $5^{34}$ $4^{36}$ $3^{42}$ $2^{20}$ $1^{27}$ $0^{44}$ $(-1)^{48}$ $(-2)^{40}$ $(-3)^{50}$ $(-4)^{52}$ $(-5)^{62}$ $(-6)^{220}$ \\
\hline
\end{tabular}
\end{center}
\vspace{-10pt}
\caption{Spectra of ${\rm SR}(4,n)$}\label{table4}
\end{table}

\medskip
The above is trivial for $m < 3$ or $n < 3$.
It was done in \cite{MW} for $m=3$, and will be done below for $n=3,4$.
The suggested spectra for $n=5$ were extrapolated from
small cases. We have not attempted to write down a proof.

\begin{Proposition}
The graph ${\rm SR}(m,3)$ has spectrum
$(3m-3)^1$, $(2m-5)^m$, $(m-3)^{m-1}$, $(m-5)^{m \choose 2}$,
$(-3)^{m(m^2-7)/6}$.
\end{Proposition}
\Proof
In view of the common part of the spectra of ${\rm SR}(m,3)$
and $J(m+2,3)$, and the fact that $m(m^2-7)/6$ is the coefficient
of $t^3$ in $\prod_{i=2}^m (1 + t + \cdots + t^{i-1})$ (for $m \ge 3$),
and the fact that the stated multiplicities sum to the total number
of vertices, it follows that we only have to show the presence of the part
$(m-3)^{m-1}$.

Fix an index $h$, $1 \le h \le m$ and consider the vector $p$
indexed by the vertices that is 1 in vertices $2e_h+e_i$
and $-1$ on vertices $e_h+2e_i$ and 0 elsewhere.
One checks that this is an eigenvector with eigenvalue $m-3$,
and the $m$ vectors defined in this way have only a single
dependency (namely, they sum to 0).
\qed

\begin{Proposition}
The graph ${\rm SR}(m,4)$ has spectrum
$(4m-4)^1$, $(3m-7)^m$, $(2m-5)^m$, $(2m-8)^{m \choose 2}$,
$(m-4)^{{m \choose 2}-1}$, $(m-6)^{m \choose 2}$,
$(m-7)^{m \choose 3}$, $(-4)^r$ where $r = {m(m^3+2m^2-13m-14)/24}$.
\end{Proposition}
\Proof
In view of the common part of the spectra of ${\rm SR}(m,4)$
and $J(m+3,4)$, and the fact that $r$ is the coefficient of $t^4$ in
$\prod_{i=2}^m (1 + t + \cdots + t^{i-1})$ (for $m \ge 4$),
and the fact that the stated multiplicities sum to the total number
of vertices, it follows that we only have to show the presence of the part
$(2m-5)^m$, $(m-4)^{{m \choose 2}-1}$, $(m-6)^{m \choose 2}$
of the spectrum.

Any eigenvector for one of these eigenvalues sums to zero on each
part of the fine equitable partition found earlier, that is, on each
set of vertices with given support. Since there are unique vertices
with support of sizes 1 or 4, these eigenvectors are 0 there,
and we need only look at the vertices $3e_i+e_j$ and $2e_i+2e_j$ and
$2e_i+e_j+e_k$.

Fix an index $h$, $1 \le h \le m$ and consider the vector $p$
(indexed by the vertices) that vanishes on each vertex
where $h$ is not in the support,
is $-1$ on $2e_h+2e_i$ and on $3e_h+e_i$,
is $2$ on $e_h+3e_i$,
is $-2$ on $2e_h+e_i+e_j$,
and is $1$ on $e_h+2e_i+e_j$.
One checks that this is an eigenvector with eigenvalue $2m-5$,
and that the $m$ vectors defined in this way are linearly independent.
That settles the part $(2m-5)^m$.

Fix a pair of indices $h,i$, $1 \le h < i \le m$,
and consider the vector $p$ (indexed by the vertices)
that is 1 on $e_h+3e_j$, $2e_i+2e_j$ and $2e_h+e_i+e_j$,
is $-1$ on $e_i+3e_j$, $2e_h+2e_j$ and $e_h+2e_i+e_j$,
and is 0 elsewhere.
One checks that this is an eigenvector with eigenvalue $m-6$,
and that the ${m \choose 2}$ vectors defined in this way
are linearly independent.
That settles the part $(m-6)^{m \choose 2}$.

Having found all desired eigenvalues except one, it is not necessary
to construct eigenvectors for the final one, since
checking $\sum \theta = {\rm tr\,} A = 0$ and
$\sum \theta^2 = {\rm tr\,} A^2 = vk$ suffices.
\qed

\end{document}